\documentclass[twocolumn,final]{autart}    

\usepackage{cite}
\usepackage{graphicx}          

\usepackage[utf8]{inputenc}
\usepackage{amssymb,amsmath,color}

\newcommand{\N}{\mathbb{N}}
\newcommand{\dist}{\mathrm{dist}} 
\renewcommand{\Re}{\mathbb{R}}
\newcommand{\smallmat}[1]{\left[ \begin{smallmatrix}#1 \end{smallmatrix} \right]}
\newcommand*\rfrac[2]{{}^{#1}\!/_{#2}}

\newcommand{\red}[1]{#1}

\graphicspath{{../figures/}}

\theoremstyle{plain}
\newtheorem{theorem}{Theorem}
\newtheorem{assumption}[theorem]{Assumption}

\theoremstyle{definition}

\theoremstyle{remark}
\newtheorem{remark}[theorem]{Remark}

\begin{document}

\begin{frontmatter}

\title{Stabilising Model Predictive Control for Discrete-time Fractional-order Systems\thanksref{footnoteinfo}} 

\thanks[footnoteinfo]{This paper was not presented at any IFAC 
meeting. Corresponding author P.~Sopasakis, Tel. +39 328 0078 230}

\author[imtl]{Pantelis~Sopasakis}\ead{p.sopasakis@imtlucca.it},\     
\author[ntua]{Haralambos~Sarimveis}\ead{hsarimv@central.ntua.gr}   

\address[imtl]{IMT Institute for Advanced Studies Lucca, 
Piazza San Fransesco 19, 55100 Lucca, Italy.}  								 
\address[ntua]{School of Chemical Engineering, National Technical 
University of Athens, 9 Heroon Polytechneiou Street, 15780 Zografou Campus, Athens, Greece.}             


\begin{abstract}                          
In this paper we propose a model predictive control scheme for constrained 
fractional-order discrete-time systems. We prove that all constraints are 
satisfied at all time instants and we prescribe conditions for the origin 
to be an asymptotically stable equilibrium point of the controlled system.
We employ a finite-dimensional approximation of the original infinite-dimensional 
dynamics for which the approximation error can become arbitrarily small. 
We use the approximate dynamics to design  a tube-based model predictive controller 
which steers the system state to a neighbourhood of the origin of controlled size. 
We finally derive stability conditions for the MPC-controlled system which are 
computationally tractable and account for the infinite dimensional nature of the 
fractional-order system and the state and input constraints.
\red{The proposed control methodology guarantees asymptotic stability of the 
discrete-time fractional order system, satisfaction of the prescribed constraints and 
recursive feasibility.}
\end{abstract}

\end{frontmatter}

\section{Introduction}
\subsection{Background and Motivation}

Derivatives and integrals of non-integer order, often referred to as
\textit{fractional}, are natural extensions of the standard
integer-order ones which enjoy certain favourable properties:
they are linear operators, preserve analyticity, and have the 
semigroup property~\cite{Pod99,Hil00}. 
Nonetheless, fractional derivatives are non-local operators, that is,
unlike integer-order ones, they cannot be evaluated
at a given point by mere knowledge of the function in 
a neighbourhood of this point and for that reason they are 
suitable for describing phenomena with infinite memory~\cite{Pod99}.

Fractional dynamics seems to be omnipresent in nature.
Examples of fractional systems include, but are not limited to, 
semi-infinite transmission 
lines with losses~\cite{ClaNarHan04}, viscoelastic 
polymers~\cite{Hil00}, magnetic core coils~\cite{Magnetism}, 
anomalous diffusion in semi-infinite bodies~\cite{GuoLiWan15} 
and biomedical applications~\cite{Magin2010} for which 
Magin \textit{et al.} provided a thorough 
review~\cite{MagOrtPodTru11}.
A good overview of the applications of fractional systems in 
physics is given in~\cite{Hil00} and~\cite{TarBook}.

A shift towards fractional-order dynamics in the field of pharmacokinetics
may be observed after the classical \textit{in-vitro-in-vivo correlations} theory 
proved to have faced its limitations~\cite{KytMacDok10}.
Non-linearities, anomalous diffusion, deep tissue trapping, 
diffusion across capillaries, synergistic and competitive 
action and other phenomena give rise to fractional-order
pharmacokinetics~\cite{DokMac08}.
In fact, Pereira derived fractional-order diffusion laws for 
media of fractal geometry~\cite{Per10}.
Increasing attention has been drawn on 
modelling and control of such systems~%
\cite{DokMac11,DokMagMac10,SopSar14}, especially in presence of
state and input constraints.

\red{%
Model predictive control (MPC) is an advanced, successful and well recognized control 
methodology and its adaptation to fractional systems is of particular interest.
The current model predictive control framework for fractional-order systems has been 
developed in a series of papers where integer-order approximations are used to formulate the 
control problem~\cite{RomDeMad+13,RomDeMad+10,BouBou10,DenCao+10,RomDeMadManVin10}. 
CARIMA (controlled auto-regressive moving average) models are often used in 
predictive control formulations for the approximation of the fractional dynamics~\cite{RomDeMad+13,RomDeMad+10,RomDeMadManVin10,Joshi2015}.
The CARIMA-based approach has been used in various applications such as the 
heating control of a semi-infinite rod~\cite{RhoBou14}, 
the power regulation of a solid oxide fuel cell~\cite{DenCao+10}
and various applications in automotive technology~\cite{RomDeMad+12}.
The celebrated Oustaloup approximation has also been used in MPC settings~\cite{RomDeMad+13}.
It should, however, be noted that such approximations aim at capturing the system
dynamics in a range of operating frequencies and, as a result, are not suitable for a
rigorous analysis and design of controllers for constrained systems.
Additionally, all of the aforementioned works provide examples of unconstrained 
systems; this shortcoming was in fact identified in the recent paper~\cite{Joshi2015}.
}

Nevertheless, this profusion of purportedly successful paradigms of MPC for 
fractional-order systems is not accompanied by a proper stability
analysis especially when input and state constraints are present.
A common denominator of all approaches in the literature is 
that they approximate the actual fractional dynamics by 
integer-order dynamics and design controllers for the 
approximate system using standard techniques. 
No stability and constraint satisfaction guarantees can be
deduced for the original fractional-order system.
Currently, one of the very few works on constrained 
control for fractional-order systems is due to Mesquine \textit{et al.}
where, however, only input constraints are taken into account \red{for 
the design of a linear feedback controller}~\cite{MesHmaBenBenTad15}.

Hitherto, two approaches can be found in the literature in 
regard to the stability analysis of discrete-time fractional
systems. The first one considers the stability of a 
finite-dimensional linear time-invariant (LTI) system, known as 
\textit{practical stability}, but fails to provide 
conditions for the actual fractional-order system to be 
(asymptotically) stable~\cite{BusKac09,GueDjeBet12}. This 
approach is tacitly pursued in many applied papers \red{where stability is
established only for a finite-dimensional approximation of the fractional-order 
system}%
~\cite{RomDeMad+13,RomTejSuar09}. 
On the other hand, fractional systems can be treated as 
infinite-dimensional systems for which various stability 
conditions can be derived (See for example~\cite[Thm.~2]{GueDjeBetMaa10}), 
but conditions are difficult to verify in 
practice let alone to use for the design of model predictive --- or other ---
controllers. 

\subsection{Contribution}
In this paper we describe a stabilising MPC framework
for fractional-order systems (of the Gr\"{u}{nwald-Letnikov type)
subject to state and input constraints.
We discretise linear continuous-time fractional dynamics
using the Gr\"{u}{nwald-Letnikov scheme which leads to infinite-dimensional linear
systems. Using a finite-dimensional approximation
we arrive at a linear time-invariant system with an additive 
uncertainty term which casts the discrepancy with the 
infinite-dimensional system. We then introduce a tube-based MPC 
control scheme which is known to steer the state to a 
neighbourhood of the origin which can become arbitrarily
small as the order of the approximation of the 
fractional-order system increases.
In our analysis, we consider both state and input constraints
which we show that are respected by the MPC-controlled system.
We finally prove that under a certain contraction-type condition
the origin is an asymptotically stable
equilibrium point for the MPC-controlled fractional-order
system (see Section~\ref{sec:stabilisation}). \red{In this work we provide,
for the first time, asymptotic stability conditions (Theorem~\ref{thm:stability-2}) and we propose a 
control methodology which guarantees the satisfaction of the prescribed state
and input constraints.}

This paper builds up on~\cite{SopNtoSar15} where the unmodelled 
part of the system dynamics was cast as a bounded
additive uncertainty term and used existing MPC theory to
drive the system's state in a neighbourhood of the origin
without, however, providing any (asymptotic) stability conditions 
for the origin.

\subsection{Mathematical preliminaries}\label{sec:mathematics}
The following definitions and notation will be used throughout
the rest of this paper.
Let $\N$, $\Re^{n}$, $\Re_+$, $\Re^{m\times n}$
denote the set of non-negative integers, 
the set of column real vectors of length 
$n$, the set of non-negative numbers and the set of 
$m$-by-$n$ real matrices respectively. 
For any nonnegative integers $k_1\leq k_2$  the finite 
set $\{k_1,\ldots,k_2\}$ is denoted by $\N_{[k_1,k_2]}$. 
Let $x$ be a sequence of real vectors of $\Re^n$.
The $k$-th vector of the sequence is denoted by
$x_k$ and its $i$-th element is denoted by $x_{k,i}$.
We denote by $\mathcal{B}_\epsilon^n=\{x\in\Re^n:\|x\|<
\epsilon\}$ the \emph{open ball} of $\Re^n$ with radius 
$\epsilon$ and we use the shorthand $\mathcal{B}^n=\mathcal{B}_1^n$. 
We define the \emph{point-to-set distance}
of a point $z\in X$ from $A$ as 
$\dist(z,A)=\inf_{a\in A}\|z-a\|$.
The space of bounded real sequences is denoted by
$\ell^\infty$.
We define the space $\ell^\infty_n$
of all sequences of real $n$-vectors $z$ so that
$(z_{k,i})_{k}\in\ell^\infty$ for $i\in\N_{[1,n]}$.

Let $\Gamma$ be a topological real vector space and $A,B\subseteq \Gamma$. 
For $\lambda\in\Re$
we define the scalar product $\lambda C=\{\lambda c: c\in C\}$
and the \emph{Miknowski sum} $A\oplus B=\{a+b: a\in A, b\in B\}$.
The Minkowski sum of a finite family of sets $\{A_i\}_{i=1}^{k}$ 
will be denoted by $\bigoplus_{i=1}^{k}A_i$. The Minkowski
sum of a sequence of sets $\{A_i\}_{i\in\N}$ is denoted by
$\bigoplus_{i\in\N}A_i$ or $\bigoplus_{i=0}^{\infty}A_i$ and is defined as the \emph{Painlev{\'e}-Kuratowski 
limit} of $\bigoplus_{i=1}^{k}A_i$ as 
$k{\to}\infty$~\cite{RocWet98}. The \emph{Pontryagin difference} between two sets 
$A,B\subseteq \Gamma$ is defined as $A\ominus B=
\{a\in A:a+b\in A, \forall b\in B\}$. A set $C$ is called
\emph{balanced} if for every $x\in C$, $-x\in C$.

\section{Fractional-order Systems}

\subsection{Discrete-time fractional-order systems}
Let $x:\Re\to\Re^n$ be a uniformly bounded function,
\textit{i.e.}, there is a $M>0$ so that $\|x(t)\|\leq M$
for all $t\in\Re$.
The Gr\"{u}nwald-Letnikov fractional-order difference of $x$ of order $\alpha>0$
and step size $h>0$ at $t$
is defined as the linear operator~\cite{OrtCoiTru15,RhoBouBou14} $\Delta^\alpha_h:\ell_n^\infty\to\ell_n^\infty$:
\begin{align}\label{eq:gl-diff}
\Delta^\alpha_h x(t) = \sum_{j=0}^{\infty}(-1)^j{\alpha \choose j}x(t-jh),
\end{align}
where ${\alpha \choose 0}=1$ and for $j\in\N$, $j>0$
\begin{align}
{\alpha \choose j}=\prod_{i=0}^{j-1}\frac{\alpha-i}{i+1}=\frac{\Gamma(\alpha+1)}{\Gamma(\alpha-j+1)j!}.
\end{align}
The \emph{forward-shifted} counterpart of $\Delta^\alpha_h$ is defined as 
${_{F}\Delta^\alpha_h} x(t)=\Delta^\alpha_h x(t+h)$.
Now, define 
\begin{align}
c^\alpha_j=(-1)^j{\alpha \choose j}= {j-\alpha-1 \choose j},
\end{align}
and notice for all $j\in\N$
that $|c^\alpha_j|\leq \alpha^j/j!$, thus, the sequence 
$(c^\alpha_j)_j$ is absolutely summable and, because of the
uniform boundedness of $x$, the series in~\eqref{eq:gl-diff}
converges, therefore, $\Delta^\alpha_h$ is well-defined.
It is worth noticing that for $\alpha\in\N$ it is $c_j^\alpha=0$
for $j\geq \lceil \alpha \rceil $, but this property does not 
hold for $\alpha \notin \N$. 
As a result, at time $t$ and for non-integer orders
$\alpha$ the whole history of $x$ is needed in order to  
estimate $\Delta^\alpha_h x(t)$.

The Gr\"{u}nwald-Letnikov difference operator gives 
rise to the Gr\"{u}nwald-Letnikov derivative 
of order $\alpha$ which is defined as~\cite[Sec.~{20}]{Sam+93}
\begin{align}
D^\alpha x(t)=\lim_{h\to 0^+}\frac{_{F}\Delta^\alpha_h x(t)}{h^\alpha}
=\lim_{h\to 0^+}\frac{\Delta^\alpha_h x(t)}{h^\alpha},
\end{align}
insofar as both limits exist. 
This derivative is then used to describe fractional-order
dynamical systems with state $x:\Re\to\Re^n$ and
input $u:\Re\to\Re^m$ as follows:
\begin{align}\label{eq:frac-sys}
\sum_{i=1}^{l} A_i D^{\alpha_i} x(t)=
 \sum_{i=1}^{r} B_i D^{\beta_i} u(t),
\end{align}
where $l,r\in\N$, $A_i$ are $B_i$ are matrices of 
opportune dimensions, all $\alpha_i$ and $\beta_i$ are nonnegative,
and by convention $D^0 x(t) = x(t)$ for any $x$.

In an Euler discretisation fashion we approximate the $D^{\alpha}$
in~\eqref{eq:frac-sys} using either $h^{-\alpha} {_{F}}\Delta^{\alpha}_h$ 
or $h^{-\alpha} \Delta^{\alpha}_h$ for a fixed step size $h$ as in~\cite{OrtCoiTru15}. 
\red{In particular, we use ${_{F}}\Delta^{\alpha}_h$ for the derivatives
of the state and $\Delta^{\alpha}_h$ for the input variables}.
We define $x_k=x(kh)$ and $u_k=u(kh)$ for $k\in\mathbb{Z}$ so the 
discretisation of~\eqref{eq:frac-sys} becomes
\begin{align}\label{eq:frac-sys-discrete}
\sum_{i=1}^{l}\bar{A}_i  \Delta_h^{\alpha_i} x_{k+1}=
\sum_{i=1}^{r}\bar{B}_i  \Delta_h^{\beta_i} u_{k},
\end{align}
with $\bar{A}_i=h^{-\alpha_i}A_i$ and
$\bar{B}_i=h^{-\beta_i}B_i$.
The involvement of infinite-dimensional operators in the
system dynamics deem these systems computationally intractable
and call for approximation methods for their simulation 
and the design of feedback controllers.

In what follows, we will approximate~\eqref{eq:frac-sys-discrete}
by a finite-dimensional state-space system treating the approximation
as a bounded additive disturbance. We then propose a control
setting which guarantees robust stability properties for%
~\eqref{eq:frac-sys-discrete}.

\subsection{Finite-dimension approximation}\label{sec:foa}
Discrete-time fractional-order dynamical
systems are essentially systems with infinite memory and an 
infinite number of state variables. As a result, standard stability theorems and 
control design methodologies for finite-dimensional systems cannot be applied directly. 
To this end we introduce the following \textit{truncated Gr\"{u}nwald-Letnikov difference operator} 
of length $\nu$ given by
\begin{align}
\Delta^\alpha_{h,\nu} x_k = 
\sum_{j=0}^{\nu}c_{j}^{\alpha}x_{k-j},
\end{align}  
System~\eqref{eq:frac-sys-discrete} is the approximated by the following
system using $\nu\geq 1$
\begin{align}\label{eq:frac-sys-approx}
{\sum_{i=1}^{l}\bar{A}_i\Delta_{h,\nu}^{\alpha_i} x_{k+1}=
\sum_{i=1}^{r}\bar{B}_i\Delta_{h,\nu}^{\beta_i} u_{k}}.
\end{align}

System~\eqref{eq:frac-sys-approx} can be written in state
space format as a linear time-invariant system with a 
proper choice of state variables $\tilde{x}_k$ as we shall explain
in this section. In the common case where the right-hand side
of~\eqref{eq:frac-sys-approx} is of the simple form 
$Bu_k$, it is straightforward to recast the system in 
state-space form. Here, we study the more general case of
equation~\eqref{eq:frac-sys-approx}, which can be written in 
the form
\begin{align}\label{eq:qwerty}
\sum_{j=0}^{\nu}\hat{A}_j x_{k-j+1} = \sum_{j=0}^{\nu}\hat{B}_j u_{k-j},
\end{align}
with $\hat{A}_j =\sum_{i=1}^{l}\bar{A}_i c_{j}^{\alpha_i}$
and $\hat{B}_j=\sum_{i=1}^{r}\bar{B}_i c_{j}^{\beta_i}$ for $j\in\N_{[0,\nu]}$.
We hereafter assume that matrix $\hat{A}_0$ is nonsingular.
With this assumption, the discrete-time dynamical system~\eqref{eq:qwerty}
becomes a \textit{normal system}, that is, future states can be determined using
past states in a unique fashion and can be written as
a linear time-invariant system~\cite[Chap.~1]{Duan10}.
Defining $\tilde{A}_j=-\hat{A}_0^{-1}\hat{A}_j$
and $\tilde{B}_j=\hat{A}_0^{-1}\hat{B}_j$,
the dynamic equation~\eqref{eq:qwerty} 
becomes
\begin{align}\label{eq:zxcv}
x_{k+1}=\sum_{j=0}^{\nu-1} \tilde{A}_j x_{k-j}
+ \sum_{j=1}^{\nu}\tilde{B}_j u_{k-j} + \tilde{B}_0 u_{k}.
\end{align}
This can be written in state space form with  state 
variable $\tilde{x}_k=(x_k, x_{k-1},\ldots, x_{k-\nu+1}, 
u_{k-1},\ldots, u_{k-\nu})'$ as
\begin{align}\label{eq:lti}
\tilde{x}_{k+1}=A\tilde{x}_{k}+Bu_{k}.
\end{align}
System~\eqref{eq:lti} is an ordinary finite-dimensional 
LTI system which will be used in the next section to 
formulate a model predictive control problem.
Throughout the rest of the paper we assume that 
the pair $(A,B)$ is stabilisable.


The truncated difference operator
$\Delta_{h,\nu}^{\alpha}$ introduces some error in the 
system dynamics. In particular, the fractional-order difference
operator ${\Delta_h^\alpha}$ can be written as
\begin{align}
\red{{\Delta_{h}^{\alpha}}=\Delta_{h,\nu}^{\alpha}+R^{\alpha}_{h,\nu}},
\end{align}
where $R^{\alpha}_{h,\nu}:\ell^{\infty}_{n}\to\ell^{\infty}_{n}$ 
is the operator $R^{\alpha}_{h,\nu}(x_k)=\sum_{j=\nu+1}^{\infty}c_{j}^{\alpha}x_{k-j}$.
Let $X$ be a compact convex subset in $\Re^n$ 
containing $0$ in its interior and at 
time $k$ assume that $x_{k-j}\in X$ for all $j\in\N$. Then, 
by the assumption that $x_{k-j}\in X$ for all $j\in\N$,
\begin{equation}\label{eq:Rnu}
R^{\alpha}_{h,\nu}(x_k)\in \bigoplus_{j=\nu+1}^{\infty}c_{j}^{\alpha}X.
\end{equation}
For all $\nu\in\N$, the right-hand side of~\eqref{eq:Rnu} is a convex
compact set with the origin in its interior.
Equation~\eqref{eq:frac-sys-discrete} can now be rewritten using the augmented
state variable $\tilde{x}$ (cf.~\eqref{eq:lti})
leading to the following linear uncertain system
\begin{align}\label{eq:uncertain-lti}
\tilde{x}_{k+1}=A\tilde{x}_{k}+Bu_{k}+Gd_k,
\end{align}
where $d_k$ is a (bounded) additive disturbance term (which depends
on $x_{k-\nu-j}$ and $u_{k-\nu-j}$ for $j\in\N$) with $G=\smallmat{I&0&\ldots&0}'$.
Assume that $u_{k-j}\in U$ for $j=1,2,\ldots$ and 
$x_{k-j}\in X$ for $j\in\N$, where $X$ and $U$ are
convex compact sets containing $0$ in their interiors.
Then, $d_k$ is bounded in a compact set $D_{\nu}$ given by
\begin{align}\label{eq:D}
D_{\nu}{=}D_{\nu}^{x}\oplus D_{\nu}^{u},
\end{align}
where
\begin{subequations}
 \begin{align} 
  D_{\nu}^{x}&=\bigoplus_{i=1}^{l} -\hat{A}_0^{-1}\bar{A}_i \bigoplus_{j=\nu+1}^{\infty} c_{j}^{\alpha_i} X,\\
  D_{\nu}^{u}&=\bigoplus_{i=1}^{r} \hat{A}_0^{-1}\bar{B}_i\bigoplus_{j=\nu+1}^{\infty} c_{j}^{\beta_i}U.
 \end{align}
\end{subequations}
Under the prescribed assumptions $D_{\nu}$ is a compact set. 
Hereafter, we shall use the notation $A^\ast_i = -\hat{A}_0^{-1}\bar{A}_i$ 
and $B^\ast_i=\hat{A}_0^{-1}\bar{B}_i$.

Recall that for a balanced set $C\subseteq\Re^n$ 
and scalars $\lambda_1, \lambda_2$ it is 
$\lambda_1 C\oplus \lambda_2 C=(|\lambda_1|+|\lambda_2|)C$.
In case $X$ and $U$ are balanced sets, the above expressions for $D_{\nu}^{x}$
and $D_{\nu}^{u}$
can be simplified. First, for $\nu\in\N$, we define the function 
$\Psi_{\nu}:\Re_+\to\Re_+$ as follows
\begin{align}
\Psi_{\nu}(\alpha)=\sum_{j=\nu+1}^{\infty}|c_j^\alpha|.
\end{align}
Then, $D_{\nu}^{x}$ is written as the finite Minkowski sum
\begin{align}\label{eq:Dn-formula}
D_{\nu}^{x} =\bigoplus_{i\in\N_{[1,l]}} A^\ast_i \Psi_{\nu}(\alpha_i) X,
\end{align}
and of course the same simplification applies to $D_{\nu}^{u}$ if $U$
is a balanced set. Notice that the computation of $D_{\nu}^{x}$
by~\eqref{eq:Dn-formula} boils down to determining a finite 
Minkowski sum,  which is possible when constraints are polytopic~\cite{GriStu93},
while overapproximations exists when they are ellipsoidal~\cite{KurEllips97}.

The size of $D_\nu$ is controlled by the choice of $\nu$;
$D_\nu$ can become arbitrarily small provided that a 
sufficiently large $\nu$ is chosen.
Notice also that $D_\nu\to \{0\}$ as $\nu\to\infty$.
In light of~\eqref{eq:uncertain-lti}, the fractional system
can be controlled by standard methods of robust control
such as min-max~\cite{DieBjo04} or tube-based MPC~%
\cite{RawMay09}; here we use the latter approach. In what follows,
we elaborate on how the tube-based MPC methodology can be
applied for the control of fractional-order systems.

\red{Various integer-order approximation methodologies have 
been proposed in the literature such as continued fraction expansions of the
system's transfer function, the approximation methods of Carlson, Matsuda,
Oustaloup, Chareff and more (see~\cite{approx2000review} for an overview).
Methods which are based on the approximation of the system dynamics in 
a given frequency range cannot lead to the formulation of an LTI system
with a bounded disturbance as in~\eqref{eq:uncertain-lti} and, as a result,
cannot be used to guarantee stability for constrained systems as we will
present in the next section.}

\section{Model Predictive Control}

\subsection{Tube-based Model Predictive Control}\label{sec:tb-mpc}
Model predictive control is a class of advanced control algorithms
where the control action is calculated at every time instant by solving
a constrained optimisation problem where a \textit{performance index}
is optimised. This performance index is used to choose an optimal sequence of 
control actions among the set of such admissible sequences, while 
corresponding state sequences are produced using a system model.
The first element of the optimal sequence is applied to the system; 
this control scheme defines the receding horizon control approach~\cite{RawMay09}.
When the process model is inaccurate, the modelling error must be taken 
into account to guarantee the satisfaction of state constraints and
closed-loop stability properties. Tube-based MPC is a flavour of MPC 
which leads to robust closed-loop stability while the accompanying 
optimisation problem is computationally tractable (unlike the min-max 
version of MPC~\cite{RawMay09}).

Here, we require that the state and input variables are 
constrained in the sets $X\subseteq \Re^n$ and 
$U\subseteq\Re^m$ respectively, both convex, compact and
contain the origin in their interior. The constraints are written 
as follows, this time involving $\tilde{x}$:
\begin{subequations}\label{eq:cstr-sys}
\begin{align}
\tilde{x}_k&\in \tilde{X},\label{eq:cstr-state}\\
u_k&\in U\label{eq:cstr-input},
\end{align}
\end{subequations}
for all $k\in\N$ and where $\tilde{X}=X^{\nu}\times U^{\nu}$, 
\textit{i.e.}, $\tilde{x}=(x_k, x_{k-1},\ldots, x_{k-\nu+1}, 
u_{k-1},\ldots, u_{k-\nu})'\in\tilde{X}$ if and only if $x_{k-i}\in X$ for $i\in\N_{[0,\nu-1]}$
and $u_{k-i}\in U$ for all $i\in \N_{[1,\nu]}$. 
Typically, in MPC $\tilde{X}$ and $U$ can be polytopes
or ellipsoids, but for our analysis no particular assumptions on $X$
and $U$ need to be imposed.


The fractional-order system is controlled
by an input $u$ which is computed according to
\begin{align}\label{eq:feedback-policy}
u_k=v_k+Ke_k,
\end{align}
where $v_k$ is a control action computed by the tube-based
MPC controller and $e_k$ is defined as the deviation between
the actual system state and the response of the nominal system,
that is $e_{k}=\tilde{x}_{k}-\tilde{z}_{k}$.
In particular, the nominal dynamics in terms of the
nominal state $\tilde{z}_k$ with input $v_k$ is
\begin{align}\label{eq:sys-nominal}
\tilde{z}_{k}=A\tilde{z}_{k-1}+Bv_{k-1}.
\end{align}

Matrix $K$ in~\eqref{eq:feedback-policy} is chosen so that
the matrix $A_K=A+BK$ is strongly stable. For $k\in \N$ let
\begin{align}
S_k^\nu=\bigoplus_{i=0}^{k} A_K^i G D_{\nu}.
\end{align}
The set $S_{\infty}^\nu=\lim_{k\to\infty} S_k^\nu$, is well-defined (the limit
exists), is compact, and is positive invariant for the deviation dynamics
$e_{k+1}=A_K e_{k}+Gd_{k}$. In what follows, $S_\infty^\nu$
will be assumed to contain the origin in its interior. For the needs
of tube-based MPC, any over-approximation of $S_{\infty}^\nu$
may be used instead~\cite{RakKerKouMay04}.

Having chosen $\tilde{z}_{0}=\tilde{x}_{0}$, it is 
$\tilde{x}_k\in \{\tilde{z}_k\}\oplus S_{\infty}^\nu$ for all
$k\in\N$. This implies that constraint~\eqref{eq:cstr-state}
is satisfied if $\tilde{z}_k\in X\ominus S_{\infty}^\nu$ and
constraint~\eqref{eq:cstr-input} is satisfied if $v_k\in U\ominus KS_{\infty}^\nu$. 
These constraints will then be involved in the 
formulation of the MPC problem which produces the control
actions $v_k=v_k(\tilde{z}_k)$.

The MPC problem amounts to the minimisation of a 
performance index $V_N$ along an horizon of future
time instants, known as the \textit{prediction horizon},
given the state of the nominal system $\tilde{z}_{k}$ at time $k$.
Let $N$ be the prediction horizon. We use the notation
$\tilde{z}_{k+i\mid k}$ for the predicted state of the nominal
system at time $k+i$ using feedback information at time $k$.
Let $\mathbf{v}_{k}=\{v_{k+i\mid k}\}_{i\in\N_{[0,N-1]}}$ 
be a sequence of input values and $\{\tilde{z}_{k+i\mid k}\}_
{i\in\N_{[1,N]}}$
the corresponding predicted states obtained by%
~\eqref{eq:sys-nominal}, \textit{i.e.}, it is
\begin{align}
\tilde{z}_{k+i+1\mid k}=A\tilde{z}_{k+i\mid k}+Bv_{k+i\mid k}, \text{ for }  i\in\N_{[0,N-1]} 
\end{align}
We introduce a performance index 
$V_N:\Re^{\bar{n}}{\times} \Re^{mN}{\to} \Re_+$ given the current state of the system
$\tilde{z}_{k\mid k}=\tilde{z}_{k}$
\begin{align}
V_N(\tilde{z}_{k\mid k}, \mathbf{v}_{k}){=}
V_f(\tilde{z}_{k+N\mid k})
{+}\sum_{i=0}^{N-1}\ell(\tilde{z}_{k+i\mid k}, v_{k+i\mid k}),
\end{align}
where $\ell$ and $V_f$ are typically quadratic functions.
We assume that $\ell(z,v)=z'Qz + v'Rv$,
where $Q$ is symmetric, positive semidefinite and $R$ is
symmetric positive definite and $V_f(z) = z'Pz$,
where $P$ is symmetric and positive definite.
The following constrained optimisation problem is then solved:
\begin{align}\label{eq:mpc_optim_prob}
\mathbb{P}_N{:}\ &V_N^\star(\tilde{z}_{k})=\min_{\mathbf{v}_{k}\in 
\mathcal{V}_N(\tilde{z}_{k})} V_N(\tilde{z}_{k}, \mathbf{v}_{k}),
\end{align}
where $\mathcal{V}_N(\tilde{z}_{k})$ is the set of all input sequences 
$\mathbf{v}_k$ with $v_{k+i\mid k}\in U\ominus KS$ for all $ i{\in}\N_{[0,N-1]}$
so that $\tilde{z}_{k+i\mid k}\in \tilde{X}\ominus S$, for all 
$i{\in}\N_{[0,N{-}1]}$ and $\tilde{z}_{k+N\mid k}\in \tilde{X}_f$ given that
$\tilde{z}_{k\mid k}=\tilde{z}_k$,
where $S$ is any over-approximation of $S_{\infty}^\nu$,
\textit{i.e.}, $S\supseteq S_{\infty}^\nu$ and 
$\tilde{X}_f\subseteq \tilde{X}$
is the \textit{terminal constraints set}.
In what follows we always assume that $\tilde{X}\ominus S$ and 
$U\ominus KS$ are nonempty sets with the origin in their
interior. In regard to the terminal cost function $V_f$ and 
the terminal constraints set $\tilde{X}_f$ we assume the following:
\begin{assumption}\label{asum:MPC_std_assumptions}
$V_f$ and $\tilde{X}_f$ satisfy the standard 
stabilising conditions in~\cite{May+00} \red{which are  
(i)  $\tilde{X}_f \subseteq \tilde{X},$ $0\in \tilde{X}_f$, $\tilde{X}_f$ is closed, 
(ii) there is a controller $\kappa_f:\tilde{X}_f\to U$ so that 
     $\tilde{X}_k$ is positively invariant for the nominal system~\eqref{eq:sys-nominal} under $\kappa_f$, i.e., 
     $A\tilde{x}+B\kappa_f(\tilde{x})\in \tilde{X}_f$ for all $\tilde{x}\in\tilde{X}_f$, and
(iii) $V_f$ is a local Lyapunov function in $\tilde{X}_f$ for the $\kappa_f$-controlled system.}
\end{assumption}

\begin{remark}\label{rmr:stability-conditions}
Matrix $P$ in $V_f$ is typically
chosen to be the (unique) solution of the 
discrete-time algebraic Riccatti equation 
$P=(A+BF)'P(A+BF)+Q+F'RF$  
with $F=-(B'PB+R)^{-1}B'PA$  and $\tilde{X}_f$ to the maximal
invariant constraint admissible set for the system $\tilde{z}_{k+1}=(A+BF)\tilde{z}_{k}$. 
Alternatively, one may choose 
$\tilde{X}_f$ to be an ellipsoid of the form $\tilde{X}_f=
\{z:V_f(z)\leq \gamma\}$ and $\gamma>0$ is chosen so that 
$\tilde{X}_f\subseteq \tilde{X}$ and $K\tilde{X}_f\subseteq U$; such 
a set can be computed according to~\cite[Sec.~{8.4.2}]{BoyVan09}.
The latter is a better choice from a computational point of
view especially in high dimensional spaces although the optimisation
problem becomes a quadratically-constrained quadratic problem.
~\hfill\ensuremath{\diamond}
\end{remark}
The solution of $\mathbb{P}_N$, namely  the optimiser 
\begin{align}
v^\star(\tilde{z}_k)=\operatorname*{argmin}\limits_{\mathbf{v}_{k}\in 
\mathcal{V}_N(\tilde{z}_{k})} V_N(\tilde{z}_{k}, \mathbf{v}_{k}),
\end{align}
defines the control law $\kappa_N(\tilde{z}_k)=v_0^\star(\tilde{z}_k)$
and leads to the closed-loop dynamics
\begin{subequations}\label{eq:composite}
\begin{align}
\tilde{x}_{k+1}&=
     A\tilde{x}_k+B\rho(\tilde{z}_k,\tilde{x}_k)+Gd_k,\\
\tilde{z}_{k+1}&=
     A\tilde{z}_{k}+B\kappa_N(\tilde{z}_{k})\label{eq:nominal}.
\end{align}
\end{subequations}
Stability properties of the closed-loop system are 
hereafter derived and stated with respect to the composite system~\eqref{eq:composite} with
state variable $(\tilde{x}, \tilde{z})$.


\subsection{Stabilising conditions}\label{sec:stabilisation}
In this section we study the stability properties of the 
controlled closed-loop system presented previously. 
Apart from the well-known stability results in robust MPC,
we prove that, under certain conditions, the controlled trajectories
of the system are asymptotically stable to the origin (see 
Theorem~\ref{thm:stability-2}). 


The following result, which readily follows 
from~\cite[Prop.~3.15]{RawMay09}, states that the system's
state converges towards $S_{\infty}$ exponentially provided
that $S=S_\infty$ is used in the formulation of the MPC
problem.
\begin{theorem}[Rawlings \& Mayne~\cite{RawMay09}]\label{thm:stability-1}
Assume that the MPC control law $\kappa_N$ stabilises
the nominal dynamical system~\eqref{eq:nominal}.
The set $S_{\infty}\times \{0\}$ is exponentially stable for
system~\eqref{eq:composite} with region of attraction
$(Z_N\oplus S_\infty)\times Z_N$, where $Z_N$ is the
domain of $\mathcal{V}_N$, \textit{i.e.}, $Z_N=\{x: 
\mathcal{V}_N(x)\neq \varnothing\}$.
\end{theorem}

In addition, the controlled trajectory of the system's
state $x_k$ and input $u_k$ satisfy constraints~%
\eqref{eq:cstr-sys} at all time instants $k\in\N$.

Notice that $S_\infty$ can become arbitrarily small with 
an appropriate choice of $\nu$ and the system's state can 
be steered this way very close to the origin, although, 
in practice large values of $\nu$ should be avoided to
limit the computation complexity of optimisation problem 
$\mathbb{P}_N$. In addition to Theorem~\ref{thm:stability-1},
we are going to prove that the state converges exactly to the 
origin and the origin is an asymptotically stable equilibrium point 
of the controlled system under certain conditions. 
The stability conditions we are about to postulate are 
easy to verify and can be used for the design of stabilising model 
predictive controllers. Hereafter, we shall assume that there are no
derivatives acting on system's inputs, \textit{i.e.}, $r=1$, $\beta_1=0$.
The main result of this section is stated as follows:

%
%
%
\begin{theorem}[Asymptotic stability]\label{thm:stability-2}
Assume that $X$ is compact and balanced, Assumption~\ref{asum:MPC_std_assumptions}
is satisfied and there is an $\epsilon\in(0,1)$ so that the following condition holds:
\begin{align}\label{eq:asym-stab-condition}
\bigoplus_{j\in\N}A_{K}^j G D&\subseteq \mathcal{B}_{\epsilon}^{\bar{n}},
\end{align}
where $D$ is the set
\begin{align}\label{eq:D2}
D= \bigoplus_{i{\in}\N_{[1,l]}}\Psi_{\nu}(\alpha_i) A^\ast_i \mathcal{B}^n.
\end{align}
Assume also that there is a $\sigma>0$ so that 
$S_\infty\subseteq \mathcal{B}_\sigma^{\bar{n}} \subseteq Z_N\oplus S_\infty$.
Then, the origin is an asymptotically stable equilibrium point for~\eqref{eq:composite}.
\end{theorem}

\begin{pf}
 The proof can be found in the appendix.
\end{pf}

%

\begin{remark}
The vector space $\Re^{\bar{n}}$ can be written as the direct sum 
of vector spaces $L_1,\ldots,L_\nu$, each of dimension $n$, so that
$\tilde{x}_k\in\Re^{\bar{n}}$ if and only if $x_{k-j+1}\in L_{j}$ for 
$j\in\N_{[1,\nu]}$.
Assume that $S_\infty\cap L_i$ has nonempty interior in the topology
of $L_i$. Then, in Theorem~\ref{thm:stability-2} one may drop the requirement that 
$S_\infty\subseteq \mathcal{B}_{\sigma}^{\bar{n}}$ by replacing the norm 
$\|\cdot\|$ of $\Re^{\bar{n}}$ by the Minkowski functional of $S_\infty$, that is 
\begin{align}
\mathrm{p}[S_\infty](\tilde{x})=\inf_{\lambda>0}\{\lambda S_\infty \ni \tilde{x}\}. 
\end{align}
The norm-ball $\mathcal{B}_\epsilon^{\bar{n}}$ becomes
$\mathcal{B}_\epsilon^{\bar{n}} = \{x: \mathrm{p}[S_\infty](x)<\epsilon\}$
and the induced matrix norm is modified accordingly, while on
$L_i$ we replace the norm by $\mathrm{p}[S_{\infty}\cap L_i](x)$.
This is based on a useful property of $\mathrm{p}[S_{\infty}\cap L_i](x)$ which 
is stated in Appendix~\ref{sec:appPropertiesPS}.~\hfill\ensuremath{\diamond}
\end{remark}

\begin{remark}
Assume that $D$ in Theorem~\ref{thm:stability-2} is a polytope (for example, 
the 1-norm or the infinity-norm is used).
Using the results presented in~\cite{RakKerKouMay04}, given a tolerance $t>0$, 
there is a $\beta>1$ and an $s\in\N$ so that the polytope 
$$F_{\beta,s}=\beta\bigoplus_{i=0}^{s}A_K^i GD$$
be an $t$-outer approximation of 
\begin{align}
 F=\bigoplus_{i=0}^{\infty}A_K^i GD,
\end{align}
in the sense that $F\subseteq F_{\beta,s} \subseteq F\oplus \mathcal{B}_t^{\bar{n}}$.
Then, the stabilising condition of Theorem~\ref{thm:stability-2} is 
satisfied if $F_{\beta,s}\subseteq \mathcal{B}_{\epsilon}^{\bar{n}}$ and this 
condition is easier to check computationally.~\hfill\ensuremath{\diamond}
\end{remark}

\begin{remark}
Since $A_K$ is a strictly Hurwitz matrix, there is a finite $a\in\N$
so that $\|A_K^{j}\|<1$ for all $j> a$. Then $F$ can be written as
\begin{align}\label{eq:x454}
 F=\bigoplus_{i=0}^{a}A_K^i GD \oplus \bigoplus_{i=0}^{\infty}A_K^{a+1+i} GD,
\end{align}
where the first term is finitely determined and in case $D$ is a polytope, it is
also a polytope. Let $\delta^\ast=\max_{d\in D}\|d\|$ (which is well-defined and finite
because $D$ is compact). The second term of $F$ in~\eqref{eq:x454} can be over-approximated 
\begin{align*}
\bigoplus_{i=0}^{\infty}A_K^{a+1+i} GD 
&\subseteq \bigoplus_{i=0}^{\infty}A_K^{a+1+i} G\mathcal{B}_{\delta^\ast}
\subseteq \bigoplus_{i=0}^{\infty} \mathcal{B}_{\delta^\ast\|A_K^{a+1}\|^i}\\
&\subseteq  \mathcal{B}_{\delta^\ast\sum_{i=0}^{\infty}\|A_K^{a+1}\|^i}
=\mathcal{B}_{\frac{\delta^\ast}{1-\|A_K^{a+1}\|}}.
\end{align*}
This is based on the observation that for a matrix 
$B\in\Re^{n\times m}$ it is $B\mathcal{B}^m\subseteq 
\mathcal{B}_{\|B\|}^n$,  where $\|B\|$ is the operator 
norm defined in Section \mbox{\ref{sec:mathematics}}. 
As a result we have that for $B_1, B_2\in\Re^{n\times m}$,
it is 
$B_1\mathcal{B}^m\oplus B_2\mathcal{B}^m\subseteq \mathcal{B}^n_{\|B_1\|}\oplus
\mathcal{B}^n_{\|B_2\|} \subseteq \mathcal{B}^n_{\|B_1\|+\|B_2\|}$.
If $a$ is adequately large and/or $\delta^\ast$ is adequately small, 
it will be $\frac{\delta^\ast}{1-\|A_K^{a+1}\|}<\epsilon<1$ (for some $\epsilon$) 
and then we can check the following stability condition
\begin{align}
  \bigoplus_{i=0}^{a}A_K^i GD \subseteq \mathcal{B}_{\epsilon-\frac{\delta^\ast}{1-\|A_K^{a+1}\|}},
\end{align}
which entails stabilising condition~\eqref{eq:asym-stab-condition} and is 
easier to verify.~\hfill\ensuremath{\diamond}
\end{remark}

\subsection{Computational complexity}
In this section we discuss the computational complexity of the proposed
scheme and give some guidelines for the selection of $\nu$. 
By Theorem~\ref{thm:stability-2}, an adequately large value of $\nu$
leads to the satisfaction of the stabilising conditions of the theorem. 
Naturally, for a given $\alpha>0$, one would be
interested to know the minimum order of approximation $\nu_\epsilon^\alpha$ for which
$\Psi_{\nu_\epsilon^\alpha}(\alpha)<\epsilon$, where $\epsilon>0$ is a desired threshold.

With $\nu=\nu_\epsilon^\alpha$, the MPC problem one needs to solve is formulated for a 
system that has $\nu_\epsilon^\alpha$ as many states as the 
original fractional system. Clearly, a parsimonious selection of $\nu$ is of
major importance for a computationally tractable controller design.
The designer needs to choose $\epsilon$ in order to strike a good
balance between performance and computational cost.
Indicatively, for $\epsilon=0.05$ and $\alpha=0.7$ we need $\nu=15$,
whereas for the same $\epsilon$ and $\alpha=1.3$ we need $\nu=4$.

\begin{figure}
\centering
\includegraphics[width=0.45\textwidth]{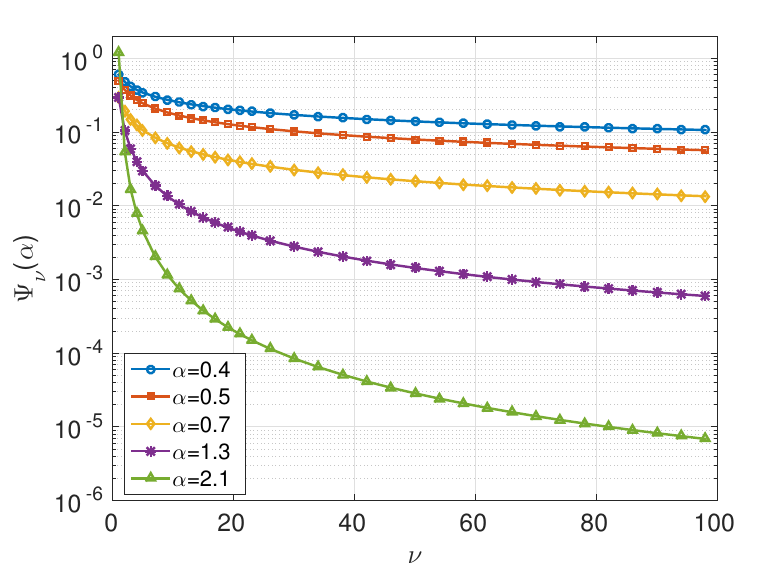}
\caption{Dependence of $\Psi_{\nu}(\alpha)$ on $\nu$ for various values of $\alpha$. 
The smaller $\alpha$ is, the slower the convergence of $\Psi_{\nu}(\alpha)$ becomes.}
\label{fig:Psinu}
\end{figure}


\section{Numerical Example}
We apply the proposed methodology to the
fractional-order system $D^\alpha x(t) = Ax(t) + Bu(t)$
with 
\begin{align}\label{eq:sys-xmpl}
A = \begin{bmatrix}
1 & 0.9\\
-0.9 & -0.2
\end{bmatrix}, B = \begin{bmatrix}
0\\
1
\end{bmatrix},
\end{align}
and $x\in\Re^2$, $u\in\Re$ and $\alpha=0.7$.
Matrix $A$ has eigenvalues $0.4\pm 0.678 i$ and the 
unactuated open-loop system is unstable.
We discretise the system with sampling period $h=0.1$
and we use $\nu=20$ based on Figure~\ref{fig:Psinu}
so that $\Psi_{\nu}(\alpha)\cong 0.041$ is adequately small (leading
to an adequately small set $D_\nu$).
This way, we derive a discrete-time LTI system
of the form $\tilde{x}_{k}=A\tilde{x}_{k-1}+Bu_{k-1}$ 
as in Section~\ref{sec:foa}.
The system state and input are subject to the constraints
\begin{subequations}\label{eq:xmpl-constr}
\begin{align}
-\begin{bmatrix}
    3\\
    3
 \end{bmatrix} \leq &x_k \leq
 \begin{bmatrix}
    3\\
    3
 \end{bmatrix},\\
-0.5 \leq &u_k \leq 0.5.
\end{align}
\end{subequations}
The terminal cost $V_f$ and the terminal constraints set 
$\tilde{X}_f$ were computed according \red{so that
Assumption~\ref{asum:MPC_std_assumptions} is satisfied}. In particular 
$\tilde{X}_f$ was chosen to be a sublevel set of $V_f$
as explained in Remark~\ref{rmr:stability-conditions}, that is
$X_f = \{x: V_f(x) \leq \gamma \}$, where $\gamma=0.015$.
The prediction horizon was chosen to be \red{$N=100$}
and the closed-loop state and input trajectories of the 
controlled system are presented in Figure~\ref{fig:sim_frac_01}
starting from the initial condition \red{$x_0=(2, 0)$}. 
Notice that the imposed constraints~\eqref{eq:xmpl-constr}
are satisfied at all time instants and the control action
saturates at its limit $u=0.5$. A phase portrait of the 
controlled system, starting from various initial points, is
shown in Figure~\ref{fig:phase_portait} and as one can see 
all trajectories converge to the origin.

In order to demonstrate the effect of $\nu$ on the system's 
closed-loop behaviour, in Figure~\ref{fig:simulations} 
we present simulations with fixed 
prediction horizon $N=100$ and different values of $\nu$
for system~\eqref{eq:sys-xmpl} starting from the initial
state~$x_0=(2, -3)$. 

\red{The average computation time for $\nu=20$ (over $150$ random (feasible) initial points $\tilde{x}_k$)
was found to be $38ms$ and the $99\%$-quantile was $43ms$ (maximum observed 
runtime: $47.3ms$). 
For a larger problem with $\nu=50$, the average runtime was
$108ms$  and the $99\%$-quantile was $195ms$ (max. $206ms$). 
The optimisation problem was formulated using 
the MATLAB toolbox YALMIP~\cite{YALMIP} and the solver MOSEK (https://www.mosek.com/).
All computations were carried out on an Intel Core i7-4510U, $4\times 2.0GHz$,
$8GB$ RAM 64-bit system running Ubuntu 14.04.
}

\begin{figure}
\centering
\includegraphics[width=0.45\textwidth]{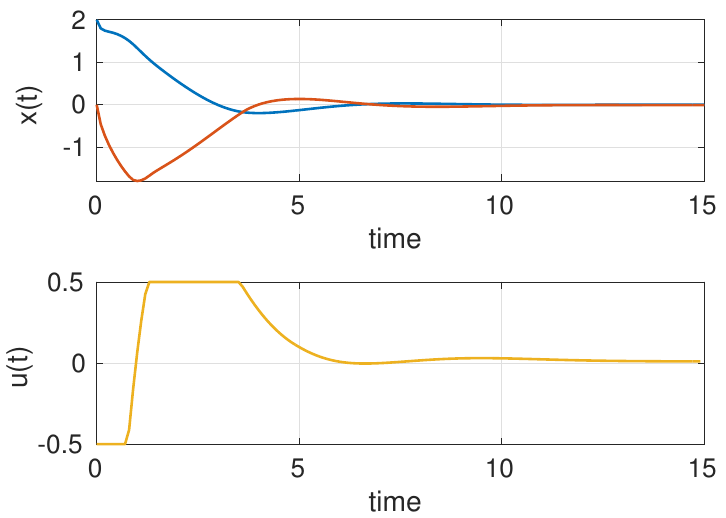}
\caption{Closed-loop simulations of system~\eqref{eq:sys-xmpl}
with the proposed MPC controller with $\nu=20$ and $N=100$. 
State (up) and input (down) trajectories.}
\label{fig:sim_frac_01}
\end{figure}

\begin{figure}
\centering
\includegraphics[width=0.45\textwidth]{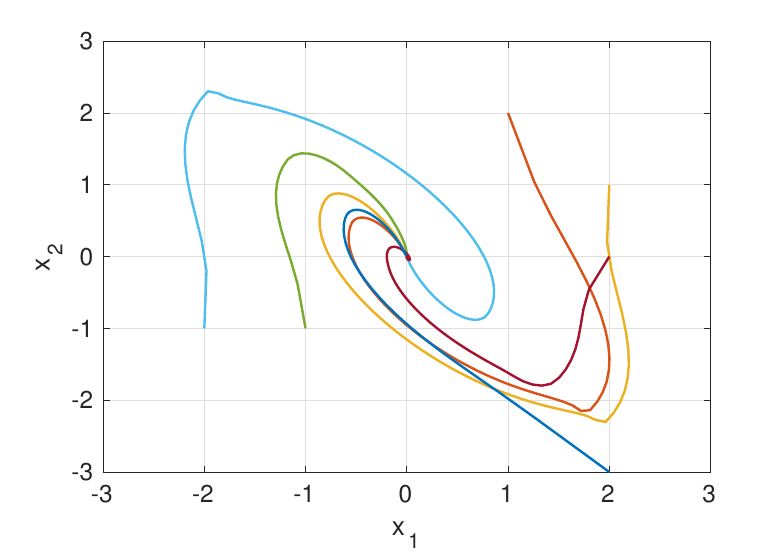}
\caption{Phase portrait of the closed-loop system starting from various different initial points \red{using $\nu=20$ and $N=100$}.}
\label{fig:phase_portait}
\end{figure}

\begin{figure}
\centering
\includegraphics[width=0.45\textwidth]{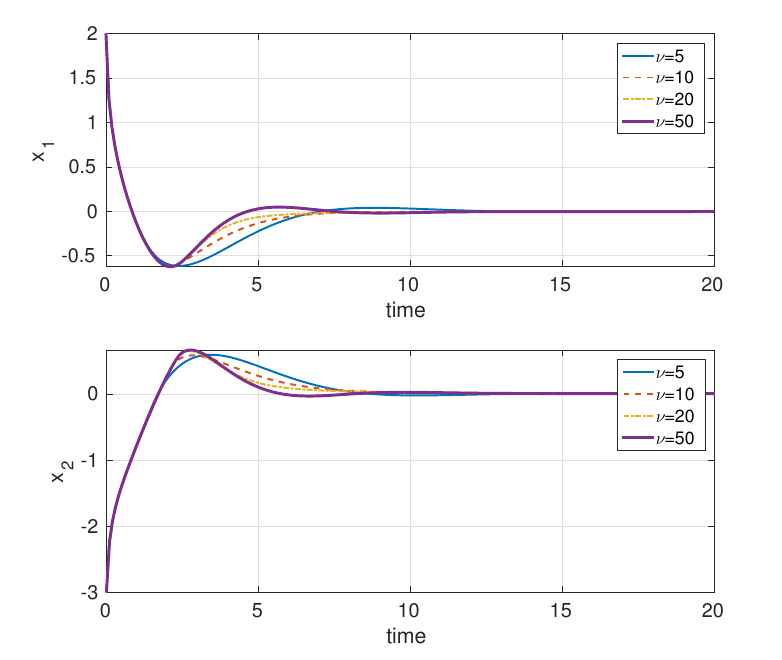}
\caption{Responses for different approximation orders $\nu$ and fixed prediction horizon $N=100$.}
\label{fig:simulations}
\end{figure}


\section{Conclusions and future work directions}
In this paper we proposed a tube-based MPC scheme for fractional 
systems which guarantees the satisfaction of state and input 
constraints. No assumptions on the fractional orders $\alpha_i$ were imposed
other than that they be nonnegative, so the results presented here are 
valid also for non-commensurate systems.
We make use of a linear and finite-dimensional approximation of the 
original dynamics and discuss how the order of approximation relates
to the computational complexity and stability properties of the resulting
controlled system.
The proposed control methodology features two important stability
properties: first, it converges exponentially fast to a convex neighbourhood
of the origin and, second, under certain conditions the origin is an 
asymptotically stable equilibrium point of the controlled system.

\red{In future work we will consider the discrepancy between the discrete-time 
fractional-order system and the original continuous-time system when the MPC control action 
is applied by a hold element. Only recently have such problems been solved for 
constrained linear time-invariant systems~\cite{SopPatSar13}.}

\begin{ack}                               
This work was funded by project 11$\Sigma\text{YN}$.10.1152, which is co-financed
by the EU and Greece, Operational Program ``Competitiveness
\& Entrepreneurship'', NSFR 2007-2013 in the context of GSRT National action ``Cooperation''.
\end{ack}

\bibliographystyle{plain}        
\bibliography{frmpc}            

\appendix
\section{Proof of Theorem~\ref{thm:stability-2}}\label{sec:appProofTHM}
We hereafter assume, without any loss of generality, that the vector-norm
$\|\cdot\|$ is the Euclidean norm and the matrix norm $\|\cdot\|$ is the 
corresponding induced norm.

(Part 1.: Attractivity)
We take $\tilde{x}_0\in Z_N$ and we shall first prove 
that the controlled trajectory of the system starting from $\tilde{x}_0$
converges to the origin (attractivity).
We start with an observation on the structure of $D_\nu$.
First, we define the function $\Phi_{\nu}(M,\alpha)=\Psi_{\nu}(\alpha)-
\Psi_{\nu+M}(\alpha)$ and notice that $D_{\nu}$ assumes
the following decomposition 
\begin{align*}
D_{\nu}=D_{\nu+M}\oplus\bigoplus_{i\in\N_{[1,l]}}\Phi_{\nu}(M,\alpha_i) A^\ast_i X,
\end{align*}
for any $M=1,2,\ldots$. Let $D^0=D_\nu$ and $S^0=S_\infty$.
Choose any $\kappa\in(\epsilon,1)$ and 
take a $0<\theta_0<\min\{1,\frac{\kappa-\epsilon}{\epsilon}\sigma\}$ and, because of 
Theorem~\ref{thm:stability-1}, there is a 
$k_0=k_0(\theta_0)\in\N$
so that for all $k\geq k_0$, $\tilde{x}_k\in S^0\oplus 
\mathcal{B}^{\bar{n}}_{\rfrac{\theta_0}{2}}\subseteq 
\mathcal{B}^{\bar{n}}_{\sigma+\theta_0}$.
Clearly, we may find $\eta_0>0$ is so that
\begin{align}\label{eq:prf1}
\bigoplus_{j\in\N}A_K^j G \mathcal{B}_{\eta_0}\subseteq \mathcal{B}_{\rfrac{\epsilon\theta_0}{2}}^{\bar{n}},
\end{align}
and take $M_0\in\N$ so that $D_{\nu+M_0}\subseteq\mathcal{B}_{\eta_0}^n$ 
and let $k\geq k_0+\nu+M_0$; then, since $\tilde{x}_k\in\mathcal{B}^{\bar{n}}_{\sigma+\rfrac{\theta_0}{2}}$, we have that 
${x}_{k-\nu-j}\in \mathcal{B}_{\sigma+\rfrac{\theta_0}{2}}^{n}$ for all $j\in\N_{[1,M_0]}$. 
Then, $d_k\in D^1$, where
$
D^1=\mathcal{B}_{\eta_0}\oplus 
\bigoplus_{i\in\N_{[1,l]}}\Phi_{\nu}(M_0,\alpha_i)A^\ast_i \mathcal{B}^n_{\sigma{+}\frac{\theta_0}{2}},
$
and refine the new target set $S^1$ as follows using the following facts (i) for all $M$ and $\nu$
it is $\Phi_\nu(M,\alpha)\leq \Psi_{\nu}(\alpha)$ (ii) condition~\eqref{eq:asym-stab-condition}
(iii) inclusion~\eqref{eq:prf1}
and (iv) because of our selection of $\theta_0$ it is $\epsilon(\theta_0+\sigma) < \kappa \sigma$.
%
%
%
%
\begin{align}
S^1&=\bigoplus_{j\in\N}A_{K}^j G D^1\\
&=\bigoplus_{j\in\N}A_{K}^j G \mathcal{B}_{\eta_0} \oplus \bigoplus_{j\in\N}A_{K}^j G 
  \bigoplus_{i\in\N_{[1,l]}}\Phi_{\nu}(M_0,\alpha_i)A^\ast_i \mathcal{B}^n_{\sigma{+}\frac{\theta_0}{2}}\notag\\
&\subseteq \mathcal{B}_{\epsilon\theta_0/2}^{\bar{n}}\oplus \bigoplus_{j\in\N}A_{K}^j G 
  \bigoplus_{i\in\N_{[1,l]}}\Psi_{\nu}(\alpha_i)A^\ast_i \mathcal{B}^n_{\sigma{+}\frac{\theta_0}{2}}\notag\\
&\subseteq \mathcal{B}_{\epsilon\frac{\theta_0}{2}}^{\bar{n}} \oplus 
  \mathcal{B}_{\epsilon(\sigma+\frac{\theta_0}{2})}^{\bar{n}} 
 \subseteq \mathcal{B}_{\epsilon(\sigma+\theta_0)}^{\bar{n}}
 \subseteq \mathcal{B}_{\kappa\sigma}^{\bar{n}}
\end{align}
and the state will converge towards $S^1$. Choose $0<\theta_1<\kappa\theta_0$. 
There is a $k_1=k_1(\theta_1)\in\N$ with
$k_1>k_0$ so that $\tilde{x}_{k}\in 
S^1\oplus \mathcal{B}_{\theta_1/2}^{\bar{n}}$ (and of course 
$\tilde{x}_k\in \mathcal{B}_{\kappa\sigma+\theta_1/2}^{\bar{n}}$)
for all $k\geq k_1$. Find $\eta_1>0$ with $\eta_1<\eta_0$ so that 
$\bigoplus_{j\in\N}A_K^j G \mathcal{B}_{\eta_1}^{\bar{n}}\subseteq 
\mathcal{B}_{\epsilon\theta_1/2}^{\bar{n}}$ and 
choose $M_1\in\N$  so that $D_{\nu+M_1}\subseteq \mathcal{B}_{\eta_1}^{n}$ 
and let $k\geq k_1+\nu+M_1$. 
Then, $x_{k-\nu-j}\in\mathcal{B}_{\kappa\sigma+\theta_1/2}^n$
for $j\in\N_{[1,M_1]}$. It follows that $d_k\in D^2$, where
$D^2=\mathcal{B}_{\eta_1}\oplus 
\bigoplus_{i\in\N_{[1,l]}}\Phi_{\nu}(M_1,\alpha_i) A^\ast_i \mathcal{B}^n_{\kappa\sigma{+}\frac{\theta_1}{2}}$
and following the same procedure as for $S^1$ we have
\begin{align}
S^2=\bigoplus_{j\in\N}A_{K}^i G D^2\subseteq 
\mathcal{B}_{\epsilon(\kappa\sigma+\theta_1)}^{\bar{n}}\subseteq 
\mathcal{B}_{\kappa^2\sigma}^{\bar{n}}
\end{align}
Recursively, we construct a sequence of sets $\{S^i\}_{i\in\N}$ 
so that $S^i\subseteq \mathcal{B}_{\kappa^i\sigma}$ and for all $i\in\N$
it is $S^i\ni 0$, therefore  $S^i\to \{0\}$ as $i\to\infty$ and, as a result, $\tilde{x}_k\to 0$
as it follows from~\cite[Ex.~{4.3(c)}]{RocWet98}.

(Part 2.: Stability)
We now need to show that the origin is stable, that is, we need to 
prove that for every $\epsilon>0$ there is a $\delta = \delta(\epsilon)>0$
so that $\|\tilde{x}_k\|<\epsilon$ for all $k\in\N$ whenever 
$\|\tilde{x}_0\|<\delta$ (and $\tilde{x}_0\in Z_N$) . 
First, notice that $\|x\|<\delta$ implies $\dist(x,S^i)<\delta$.
By Theorem~\ref{thm:stability-1} we know that for each $i\in\N$ and given $\epsilon$ there is a 
$\delta^\ast={\delta}^\ast(\epsilon, i)$ so that $\dist(x_0,S^i)<{\delta}^\ast$ implies
$\dist(x_k,S^i)<\rfrac{\epsilon}{2}$ for all $k\in\N$. 

Let $i=i(\epsilon)=\lceil \log_{\kappa}\frac{\epsilon}{2\sigma}\rceil$ and
take $x_0$ so that $\|x_0\| < \delta(\epsilon, i(\epsilon))$; then 
$\dist(x_0, S^i)<\delta(\epsilon, i(\epsilon))$, therefore 
for all $k\in\N$, $\dist(x_k, S^i)<\frac{\epsilon}{2}$ 
and $\|x_k\|<\frac{\epsilon}{2}+\kappa^i \sigma<\epsilon$.
~\hfill\ensuremath{\square}

\section{Properties of $\mathrm{p}[S]$ and $\mathrm{p}[S\cap L_i]$}\label{sec:appPropertiesPS}
Let $\beta=\max_{s\in S}\|s\|$. Then $S\subseteq \mathcal{B}_{\beta}$,
thus for all $\tilde{x}\in\Re^{\bar{n}}$, $\mathrm{p}[S](x)\geq \mathrm{p}[\mathcal{B}_{\beta}^{\bar{n}}](\tilde{x})$,
\textit{i.e.}, $\mathrm{p}[S](\tilde{x})\geq \frac{\|\tilde{x}\|}{\beta}$, or equivalently $\|\tilde{x}\|\leq \beta \mathrm{p}[S](\tilde{x})$.
Given that $S\cap L_i$ has nonempty interior, we may find $\gamma_i>0$ with $\mathcal{B}_{\gamma_i}^n \subseteq S\cap L_i$. 
Let $x$ be the projection of $\tilde{x}$ on $L_i$. We then have $\mathcal{B}^n_{\gamma_i}\subseteq S \cap L_i$, thus
$\mathrm{p}[S \cap L_i](x) \leq \mathrm{p}[\mathcal{B}_{\gamma_i}^n](x) =  \frac{\|x\|}{\gamma_i}$.
We then have
$
\gamma_i \mathrm{p}[S \cap L_i](x) \leq \|x\| \leq \|\tilde{x}\| \leq \beta \mathrm{p}[S](\tilde{x})
$
therefore
$
\mathrm{p}[S \cap L_i](x) \leq a_i \mathrm{p}[S](\tilde{x}),
$
with $a_i=\frac{\beta}{\gamma_i}$.~\hfill\ensuremath{\square}
\end{document}